\documentclass[12pt,leqno,a4paper,final]{amsart}
\usepackage{amsmath, amsthm, verbatim, amsfonts, amssymb, color}
\usepackage{graphicx,epic,eepic}

\usepackage{mathrsfs,graphicx}
\usepackage{dsfont}
\usepackage[a4paper,margin=1in]{geometry}
\usepackage{marginnote}
\usepackage[notcite,notref]{showkeys}

\theoremstyle{plain}
\newtheorem*{thm}{Theorem}

\theoremstyle{definition}
\newtheorem{remark}{Remark}

\theoremstyle{remark}

\newcommand\Hup{{\mathds{H}^\uparrow}}
\newcommand\nat{\mathds N}
\newcommand\real{\mathds R}
\newcommand\comp{\mathds C}
\newcommand\I{\mathds 1}

\newcommand\RE{\operatorname{Re}}
\newcommand\IM{\operatorname{Im}}

\newcommand{\et}{\quad\text{and}\quad}

\newcommand\dup{\mathrm{d}}
\newcommand\eup{\mathrm{e}}

\begin{document}\allowdisplaybreaks
\title[Complete Bernstein functions \& subordinators with nested ranges]{Complete Bernstein functions and subordinators with nested ranges.\\A note on a paper by P.\ Marchal}

\author[C.-S.~Deng]{Chang-Song Deng}
\address[C.-S.~Deng]{School of Mathematics and Statistics\\ Wuhan University\\ Wuhan 430072, China}
\email{dengcs@whu.edu.cn}
\thanks{The first-named author gratefully acknowledges support through the Alexander-von-Humboldt foundation, the National Natural Science Foundation of China (11401442) and the China Postdoctoral Science Foundation (2016T90719).}

\author[R.\,L.~Schilling]{Ren\'e L.\ Schilling}
\address[R.\,L.~Schilling]{TU Dresden\\ Fachrichtung Mathematik\\ Institut f\"{u}r Mathematische Stochastik\\ 01062 Dresden, Germany}
\email{rene.schilling@tu-dresden.de}


\date{\today}

\begin{abstract}
    Let $\alpha:[0,1]\to [0,1]$ be a measurable function. It was proved by P.\ Marchal \cite{Mar15} that the function
    $$
        \phi^{(\alpha)}(\lambda):=\exp\left[
        \int_0^1\frac{\lambda-1}{1+(\lambda-1)x}\,\alpha(x)\,\dup x
        \right],\quad \lambda>0
    $$
    is a special Bernstein function. Marchal used this to construct, on a single probability space, a family of regenerative sets $\mathcal R^{(\alpha)}$ such that $\mathcal{R}^{(\alpha)} \stackrel{\text{law}}{=} \overline{\{S^{(\alpha)}_t:t\geq 0\}}$ ($S^{(\alpha)}$ is the subordinator with Laplace exponent $\phi^{(\alpha)}$) and $\mathcal R^{(\alpha)}\subset \mathcal R^{(\beta)}$ whenever $\alpha\leq\beta$.
    We give two simple proofs showing that $\phi^{(\alpha)}$ is a complete Bernstein function and extend Marchal's construction to all complete Bernstein functions.
\end{abstract}

\maketitle

For a measurable function $\alpha:[0,1]\to [0,1]$ define
\begin{equation}\label{function}
    \phi^{(\alpha)}(\lambda):=\exp\left[
        \int_0^1\frac{\lambda-1}{1+(\lambda-1)x}\,\alpha(x)\,\dup x
        \right],\quad\lambda>0.
\end{equation}
If $\alpha\equiv\alpha_0$ is a constant function with $\alpha_0\in(0,1)$, then $\phi^{(\alpha)}$ reduces to the fractional power function $\lambda\mapsto\lambda^{\alpha_0}$. In a recent paper, Marchal \cite{Mar15} proves that for any measurable function $\alpha:[0,1]\to [0,1]$, $\phi^{(\alpha)}$ is
a special Bernstein function, and the dual Bernstein function $\lambda/\phi^{(\alpha)}(\lambda)$ is $\phi^{(1-\alpha)}$. As an application, Marchal constructs, on a single probability space, a family of regenerative sets $\mathcal R^{(\alpha)}$ such that $\mathcal{R}^{(\alpha)} \stackrel{\text{law}}{=} \overline{\{S^{(\alpha)}_t:t\geq 0\}}$ ($S^{(\alpha)}$ is the subordinator with Laplace exponent $\phi^{(\alpha)}$) and $\mathcal R^{(\alpha)}\subset \mathcal R^{(\beta)}$ whenever $\alpha\leq\beta$.
In this short note, we will go further to show that $\phi^{(\alpha)}$ is a complete Bernstein function for all measurable weights $\alpha:[0,1]\to [0,1]$ and that Marchal's construction holds for all complete Bernstein
functions. Independently of us this has
been remarked by Alili, Jedidi and Rivero
in \cite[Example 4.2, p.~730]{AJR14}.

Let us first briefly recall some basic facts on Bernstein functions. We use the monograph \cite{SSV12} as our standard reference for Bernstein functions. A function $f:(0,\infty)\to [0,\infty)$ is called a Bernstein function, if $f\in C^\infty(0,\infty)$ and $(-1)^{k-1}f^{(k)}\geq0$ for all $k\in\nat$. All Bernstein
functions admit a unique L\'evy--Khintchine representation
\begin{equation}\label{repr}
    f(\lambda)=a+b\lambda+\int_{(0,\infty)}
    \left(1-\eup^{-\lambda x}\right)\,\nu(\dup x),
\end{equation}
where $a,b\geq0$ and $\nu$ is a 
measure on $(0,\infty)$ satisfying $\int_{(0,\infty)}(x\wedge1)\,\nu(\dup x)<\infty$. A Bernstein function $f$ is said to be a special Bernstein function if $f^*(\lambda):=\lambda/f(\lambda)$ is again a Bernstein function; in this case, $f^*$ is called the dual Bernstein function of $f$. A Bernstein function $f$ is a complete Bernstein function if its L\'{e}vy measure $\nu$ in \eqref{repr} has a completely monotone density $m$ (i.e.\ $m\in C^\infty(0,\infty)$ and $(-1)^{k}m^{(k)}\geq0$ for all $k\in\nat\cup\{0\}$) w.r.t.\ Lebesgue measure. We use $\mathcal{BF}$, $\mathcal{SBF}$ and $\mathcal{CBF}$ to denote the collections of all Bernstein functions, special Bernstein functions and complete Bernstein functions, respectively. It is known that
$$
    \mathcal{CBF}
    \subsetneqq\mathcal{SBF}
    \subsetneqq\mathcal{BF},
$$
see \cite[Propositions 11.16 and 11.17 and Example 11.18]{SSV12}. In  contrast to $\mathcal{SBF}$, the class $\mathcal{CBF}$ has well-understood structural properties and many examples of complete Bernstein functions are
known, cf.\ \cite[Chapters 6 and 16]{SSV12}.

We can now state the main result of this note.
\begin{thm}
    For any measurable function $\alpha:[0,1]\to [0,1]$, the function $\phi^{(\alpha)}$ defined by \eqref{function} is a complete Bernstein function.
\end{thm}
\begin{remark}\label{rem-1}
    Let $c,d\in\real$ with $c<d$. For any measurable function $\alpha:[c,d]\to [0,1]$, it follows from our theorem and
    a straightforward change of variables that the function
    $$
        \lambda\mapsto\exp\left[
        \int_c^d\frac{\lambda-1}{(d-c)+(\lambda-1)(x-c)}
        \,\alpha(x)\,\dup x
        \right],\quad \lambda>0
    $$
    is also a complete Bernstein function.
\end{remark}

\begin{remark}\label{rem-2}
    Our second proof of the theorem shows, in particular, that -- up to a multiplicative constant $c>0$ -- \emph{\textbf{all} complete Bernstein functions have a representation of the form \eqref{function}; moreover, the function $\alpha(x)$ is uniquely determined by the corresponding Bernstein function
    and vice versa}.

    For this we use the following characterization of complete Bernstein functions, see \cite[Theorem 6.17]{SSV12}. We have $f\in\mathcal{CBF}$ if, and only if, there is some $\gamma\in\real$ and a measurable function $\eta:[0,\infty)\to[0,1]$
    such that
    \begin{equation}\label{comp-rep}
        f(\lambda)
        = \exp\left[\gamma+\int_0^\infty \left(\frac{t}{1+t^2}-\frac{1}{\lambda+t}\right)\eta(t)\,\dup t\right],\quad \lambda>0.
    \end{equation}
    The pair $(\gamma,\eta)$ uniquely characterizes $f\in\mathcal{CBF}$ and vice versa.

    We will see that there is a one-to-one correspondence $\eta\leftrightarrow\alpha$ given by $\eta(t) = \alpha\left(\frac 1{1+t}\right)$, $t\in [0,\infty)$, while $\gamma=\gamma(\alpha)$. Since $\phi^{(\alpha)}(1)=1$, this means that any $f\in\mathcal{CBF}$ can be written as $f(1)\times\phi^{(\alpha)}$. At the level of subordinators this amounts to consider the time-changed subordinator $(S^{(\alpha)}_{ct})_{t\geq 0}$, $c=f(1)>0$; obviously, $\{S^{(\alpha)}_{ct}:t\geq 0\}=\{S^{(\alpha)}_t:t\geq 0\}$, i.e.\ Marchal's Theorem 2 holds for \emph{all} complete Bernstein functions.
\end{remark}

\section{First Proof}
Our first proof relies on the fact that $f\in\mathcal{CBF}$ if, and only if, $f$ has an analytic extension onto the open upper complex half-plane $\Hup:=\{z\in\comp:\IM z>0\}$ such that $f:\Hup\to\Hup$ and and $f(0+)=\lim_{(0,\infty)\ni\lambda\to 0}f(\lambda)$ exists, see~\cite[Theorem~6.2]{SSV12}.

\begin{proof}[First proof of the main theorem]
    According to \cite[Theorem I.4.3, p.~32]{SS05},
    we can pick a sequence of step functions $\{\alpha_n:n\in\nat\}$ on $[0,1]$ of the following form
    $$
        \alpha_n=\sum_{i=1}^na^{(n)}_{i} \I_{\left[t^{(n)}_{i-1},t^{(n)}_{i}\right)},
    $$
    where $a^{(n)}_{i}\in[0,1]$ for all $i\in\{1,\dots,n\}$ and $0=t^{(n)}_0<t^{(n)}_1<\dots<t^{(n)}_n=1$, such that $\alpha_n(x)\to \alpha(x)$ for almost all $x\in[0,1]$ as $n\to \infty$. For $n\in\nat$ and $\lambda>0$, we have
    \begin{align*}
        \phi^{(\alpha)}_n(\lambda)
        :=& \exp\left[  \int_0^1\frac{\lambda-1}{1+(\lambda-1)x}\,\alpha_n(x)\,\dup x \right]\\
        =& \exp\left[ \sum_{i=1}^n a^{(n)}_{i}
        \int_{t^{(n)}_{i-1}}^{t^{(n)}_{i}}
        \frac{\lambda-1}{1+(\lambda-1)x}\,\dup x\right]\\
        =& \prod_{i=1}^n\left( \frac{1+(\lambda-1)t^{(n)}_{i}}{1+(\lambda-1)t^{(n)}_{i-1}}\right)^{a^{(n)}_{i}}.
    \end{align*}
    This representation allows us to extend $\phi^{(\alpha)}_n$ analytically onto the open upper half-plane $\Hup$.
    Moreover,
    $$
        \lim_{(0,\infty)\ni\lambda\to 0}
        \phi^{(\alpha)}_n(\lambda)=0\quad
        \text{for all $n\in\nat$},
    $$
    and by the dominated convergence theorem, one has
    $$
        \lim_{n\to \infty}\phi^{(\alpha)}_n(\lambda)
        =\phi^{(\alpha)}(\lambda)\quad
        \text{for all $\lambda>0$}.
    $$
    Let $n\in\nat$ and $z\in\Hup$. Note that
    \begin{align*}
        \phi^{(\alpha)}_n(z)
        &= \prod_{i=1}^n\left( \frac{1+(z-1)t^{(n)}_{i}}{1+(z-1)t^{(n)}_{i-1}} \right)^{a^{(n)}_{i}}\\
        &= \left( t^{(n)}_{1} \right)^{a^{(n)}_{1}}
           \left( z-1+\frac{1}{t^{(n)}_{1}} \right)^{a^{(n)}_{1}}
            \prod_{i=2}^n \left( \frac{t^{(n)}_{i}}{t^{(n)}_{i-1}} \right)^{a^{(n)}_{i}}
            \left( \frac{z-1+{1}/{t^{(n)}_{i}}}{z-1+{1}/{t^{(n)}_{i-1}}} \right)^{a^{(n)}_{i}}\\
        &= \left( t^{(n)}_{1} \right)^{a^{(n)}_{1}}
           \left( \frac{z^{(n)}_{1}}{z^{(n)}_{0}} \right)^{a^{(n)}_{1}}
            \prod_{i=2}^n\left( \frac{t^{(n)}_{i}}{t^{(n)}_{i-1}} \right)^{a^{(n)}_{i}}
            \left( \frac{z^{(n)}_{i}}{z^{(n)}_{i-1}} \right)^{a^{(n)}_{i}},
    \end{align*}
    where
    $$
        z^{(n)}_{0}:=1,\quad
        z^{(n)}_{i}:=z-1+\frac{1}{t^{(n)}_{i}},\quad i=1,\dots,n.
    $$
    It is easy to see that $\IM z^{(n)}_{i}=\IM z$ for
    all $i\in\{1,\dots,n\}$, and $\RE z^{(n)}_{i-1}>\RE z^{(n)}_{i}$
    for all $i\in\{2,\dots,n\}$; see the figure below.
    \begin{figure}[ht]
        \includegraphics[width = .7\textwidth]{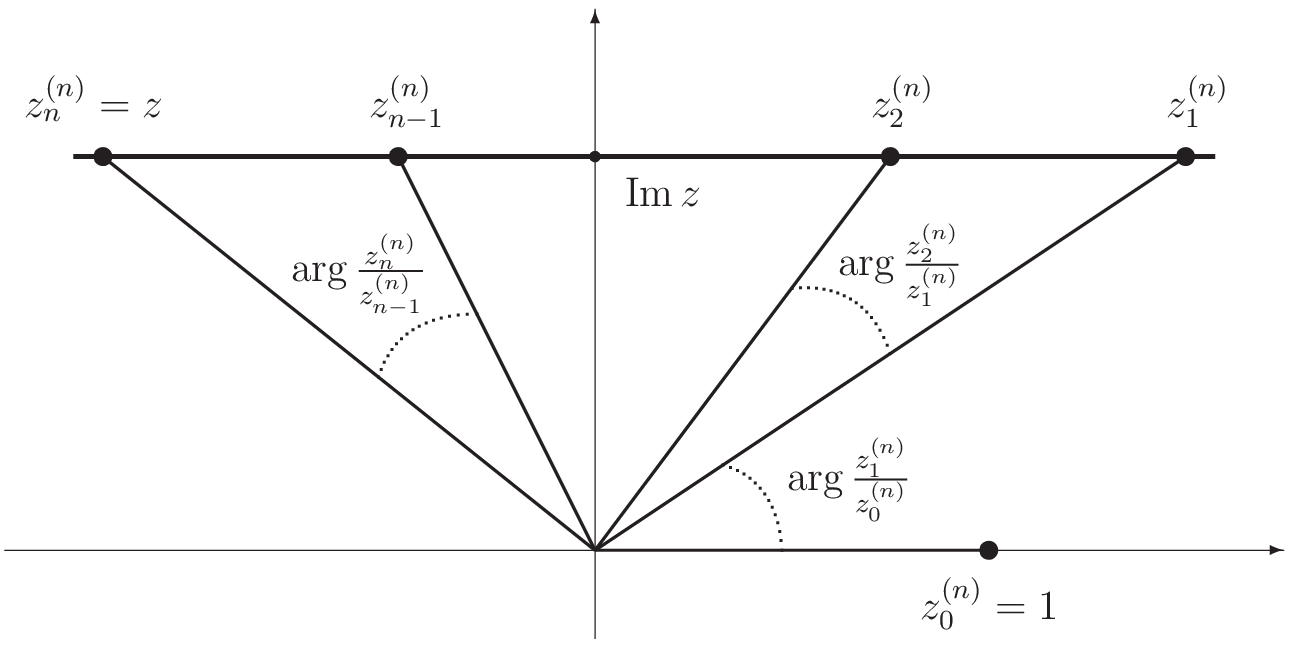}
    \end{figure}

\noindent Then we find that
    \begin{equation}\label{ds3s}
        \arg\frac{z^{(n)}_{i}}{z^{(n)}_{i-1}}>0
        \;\;\text{for every $i\in\{1,\dots,n\}$,}\;\; \text{and}
        \;\;
        \sum_{i=1}^n\arg\frac{z^{(n)}_{i}}{z^{(n)}_{i-1}}
        = \arg z^{(n)}_{n}
        = \arg z\in(0,\pi).
    \end{equation}
    If $a^{(n)}_{i}=0$ for all $i\in\{1,\dots,n\}$, then
    $\phi^{(\alpha)}_n\equiv1$, which is obviously
    of class $\mathcal{CBF}$; otherwise, we obtain
    from \eqref{ds3s} and $a^{(n)}_{i}\in[0,1]$ that
    $$
        \arg \phi^{(\alpha)}_n(z)
        =\sum_{i=1}^na^{(n)}_{i}\arg\frac{z^{(n)}_{i}}{z^{(n)}_{i-1}}
        \in (0,\pi),
    $$
    which implies that $\phi^{(\alpha)}_n$ preserves the open
    upper half-plane, and
    so $\phi^{(\alpha)}_n\in\mathcal{CBF}$. Thus,
    we conclude that $\phi^{(\alpha)}_n\in\mathcal{CBF}$ for all $n\in\nat$. Since $\mathcal{CBF}$ is closed under pointwise limits, cf.\ \cite[Corollary 7.6\,(ii)]{SSV12}, it follows that $\phi^{(\alpha)}\in\mathcal{CBF}$.
\end{proof}

\section{Second Proof}

Our second proof is based on the characterization \eqref{comp-rep} of complete Bernstein functions mentioned in Remark~\ref{rem-2}. Alili, Jedidi and Rivero have discovered the same argument, independently of us in \cite[Example 4.2, p.~730]{AJR14}. Since their proof appears in a different context and contains a small mistake, we provide the short proof for the readers' convenience. We are grateful to an anonymous referee pointing out the reference \cite{AJR14} and we acknowledge their priority for this argument.

\begin{proof}[Second proof of the main theorem]
    Observe that
    \begin{align*}
        \log\phi^{(\alpha)}(\lambda)&=
        \int_0^1 \frac{\lambda-1}{1+(\lambda-1)x}\,\alpha(x)\,\dup x\\
        &=
        \int_0^1 \frac 1{x^2} \left(x -
        \frac{1}{(1-x)/x+ \lambda}\right)\alpha(x)\,\dup x,\quad \lambda>0.
    \end{align*}
    Changing variables according to $t =(1-x)/x$ yields
    that for $\lambda>0$
    \begin{align*}
        \log\phi^{(\alpha)}(\lambda)
        &= \int_0^\infty \left(\frac 1{1+t} - \frac{1}{\lambda+t }\right)\alpha\left(\frac 1{1+t}\right)\,\dup t\\
        &=\int_0^\infty \left(\frac 1{1+t}-\frac{t}{1+t^2 }\right)\alpha\left(\frac 1{1+t}\right)\,\dup t\\
        &\quad+\int_0^\infty \left(\frac t{1+t^2} - \frac{1}{\lambda+t }\right)\alpha\left(\frac 1{1+t}\right)\,\dup t.
    \end{align*}
    Since
    \begin{align*}
    \left|\frac t{1+t^2} - \frac{1}{\lambda+t}\right|
    = \left|\frac{\lambda t-1}{\lambda+t}\right| \frac{1}{1+t^2} \in L^1((0,\infty);\dup t)\quad \text{for all $\lambda>0$},
    \end{align*}
    we know that both integrals appearing in the above
    representation of $\log\phi^{(\alpha)}(\lambda)$ are finite.
    This shows that $\phi^{(\alpha)}$ is a complete Bernstein function of the form \eqref{comp-rep} with parameters
    \begin{gather*}
        \gamma := \int_0^\infty \left(\frac 1{1+t}-\frac{t}{1+t^2 }\right)\alpha\left(\frac 1{1+t}\right)\,\dup t
        \et
        t\mapsto \eta(t) := \alpha\left(\frac 1{1+t}\right).
    \qedhere
    \end{gather*}
\end{proof}

\end{document}